\input amstex
\input amsppt.sty
\loadbold
\magnification=1200
\def\la{{\lambda}}
\def\bla{{\boldsymbol \lambda}}

\def\sD{{\Cal D}}

\def\sH{{\Cal H}}

\def\sR{{\Cal R}}

\def\sY{{\Cal Y}}

\def\bH{{\bold H}}

\def\bT{{\bold T}}

\def\BZ{{\Bbb Z}} 

\def\ph{\varphi}

\def\la{\lambda}

\def\La{\Lambda}

\def\le{\leqslant}
\def\ge{\geqslant}
\def\fS{\frak S}\def\fT{\frak T}

\def\res{\text{\rm res}}

\def\Hom{\text{\rm Hom}}

\def\lecq{\preccurlyeq}

\def\recq{\succcurlyeq}

\def\bft{{\bold t}}

\def\fku{{\frak u}}
\def\bka{{\boldkey a}}
\def\bkb{{\boldkey b}}
\def\bkc{{\boldkey c}}
\def\bkM{{\boldkey M}}

\def\bla{{\boldsymbol \lambda}}
\def\brho{{\boldsymbol \rho}}
\def\bmu{{\boldsymbol \mu}}
\def\bka{{\boldkey a}}
\def\bkb{{\boldkey b}}
\def\bkc{{\boldkey c}}

\def\rn{{\frak n}}
\def\proclaim#1{\vskip.3cm\noindent{\bf#1.}\quad\it}
\def\endproclaim{\par\vskip.3cm\rm}
\def\tpi{{\tilde\pi}}

\topmatter
\title Specht Modules and Branching Rules for Ariki-Koike Algebras\endtitle
\rightheadtext{Ariki-Koike  algebras}
\thanks
Research was supported by the 
Australian Research Council. The second author was also supported by 
Sonderforschungsbereich 343 at the University of Bielefeld and the
National Natural Science Foundation in China. He
wishes to thank the University of Bielefeld for its hospitality during
his visit.
\endthanks
\author Jie Du  and Hebing Rui  \endauthor
\affil {\eightpoint
School of Mathematics, University of New South Wales\\  Sydney,
 2052, Australia \\ 
 Department of Mathematics, University of Shanghai for Science \& Technology,\\
Shanghai, 200093, China}\endaffil 
\date 18 January, 1999\enddate
\email  jied\@maths.unsw.edu.au and hbruik\@online.sh.cn\endemail
\endtopmatter
\vsize24truecm

Specht modules for an Ariki-Koike algebra $\bH^{r}_m $
have been investigated recently
in the context of cellular algebras (see, e.g., \cite{GL} and \cite{DJM}). 
Thus, these modules 
are  defined as quotient modules of  certain ``permutation'' modules, that is,
defined as ``cell modules'' via cellular bases. So 
cellular bases play a decisive r\^ole in these work.
However, the classical theory \cite{C} or the work in the case
when $m=1,2$ (i.e., the case for type $A$ and $B$ Hecke algebras)
suggest that a construction as {\it submodules} 
without using cellular bases should exist.
Following our previous work \cite{DR}, we shall introduce in this paper
Specht modules for an Ariki-Koike algebra as submodules of those
``permutation'' modules and investigate their basic properties
such as Standard Basis Theorem and the ordinary Branching Theorem,
generalizing several classical constructions given in 
\cite{JK} and \cite{DJ} for type $A$.

The second part of the paper moves on looking for modular branching rules 
for Specht and irreducible modules over an Ariki-Koike algebra.
These rules for symmetric groups were recently established
by Kleshchev \cite{K}, and were generalized to Hecke algebras of type $A$ by
Brundan \cite{B}. We shall restrict to the case where $\bH^{r}_m $ has a 
semi-simple bottom in the sense of \cite{DR}. This is because
the classification of irreducible modules is known in this case.
Our Specht module
theory and the Morita equivalence theorem established in \cite{DR}
are the main ingredients in this generalization.

We point out that this realization as submodules is actually very important 
in the Specht/Young module theory for Ariki-Koike algebras and their
associated endomorphism algebras. See \cite{CPS} for more details.
We emphasise that 
our method throughout the work is independent of the use of cellular bases.
Moreover, in the proofs of the main results (2.2), (3.6), (4.2) and (4.10), one
will see how the relevant level 1 results (i.e., the results for type $A$
Hecke algebras)
are ``lifted'' to the corresponding level $m$ results.

After a first manuscript of the paper was completed, the authors 
received a preprint entitled ``On the classification
of simple modules for cyclotomic Hecke algebras of type $G(m,1,n)$
and Kleshchev multipartitions'' by S. Ariki. With this latest development,
one would expect the branching rules (4.10) holds in general.

\subhead 1. Specht modules\endsubhead
Let $\fS_r=\fS_{\{1,\cdots,r\}}$ be the symmetric group
of all permutations of $1,\cdots,r$ with Coxeter generators  $s_i=(i,i+1)$,
 and let $\fS_\la$ denote the Young
subgroup corresponding the composition $\la$ of $r$.
(A {\it composition} $\la$ of $r$  is a finite sequence 
  of non-negative integers $(\la_1, \la_2, \cdots,\la_m )$
such that  $|\la|=\sum_i\la_i=r$.) Thus,  we have
$$
\fS_\la=\fS_\bka=\frak S_{\{1,\cdots,a_1\}}\times 
\frak S_{\{a_1+1, \cdots, a_2\}}
\times \cdots \times  \frak S_{\{a_{m-1}+1, \cdots, r\}},\tag1.1
$$
where $\bka=[a_0,a_1,\cdots,a_m]$ with $a_0=0$ and
$a_i=\la_1+\cdots+\la_i$ for all $i=1,\cdots,m$.
We will denote by $\sD_\la$ the set of distinguished representatives of right
$\frak S_\la$-cosets and write $\sD_{\la\mu}=\sD_\la\cap\sD_\mu^{-1}$, which
is
the set of distinguished representatives of double cosets
$\frak S_\la\backslash\fS_r/\fS_\mu$.

For later use, let $\La(r)$ (resp. $\La^+(r)$) denote the set of all 
compositions (resp. partitions) of $r$. For $\la\in\La^+(r)$, let
$\la'$ be the dual partition of $\la$: thus $\la'_i=\#\{\la_j\ge i\}$.
There is a unique element $w_\la\in\sD_{\la\la'}$ with the trivial
intersection property: 
$$\fS_\la^{w_\la}\cap \fS_{\la'}=w_\la^{-1}\fS_\la w_\la\cap \fS_{\la'}=\{1\}.\tag1.2$$
If we represent a partition by a Young diagram $\sY(\la)$, for example,
 $\sY(\la)=\boxed{\smallmatrix \square&\square&\square\cr\square&\square&\cr
\endsmallmatrix}$ represents $\la=(32)$, then $w_\la$
is defined by $\bft^\la w_\la=\bft_\la$, where
$\bft^\la$ (resp. $\bft_\la$) is the $\la$-tableau obtained
by putting the numbers $1,2,\cdots,r$ in order 
into the boxes from left to right
down successive rows (resp. columns): thus
$\bft^{(32)}=\boxed
{\smallmatrix1 &2&3\cr4&5&\cr
\endsmallmatrix}$, and
$\bft_{(32)}=\boxed{
\smallmatrix 1&3&5\cr2&4&\cr\endsmallmatrix}$.

An $m$-tuple $\bla=(\la^{(1)},\cdots,\la^{(m)})$ of 
 partitions with $r=\sum_i|\la^{(i)}|$ is called a multi-partition
or an
$m$-partition of $r$. Define $m$-composition similarly.
Let $\bla'=(\la^{(m)\prime},\cdots,\la^{(1)\prime})$ denote the
$m$-composition {\it dual} to $\bla$.
By concatenating the components of $\bla$, the resulting composition
of $r$ will be denoted by 
$$\bar\bla=\la^{(1)}\vee\cdots\vee\la^{(m)}.$$ Thus,
the bar $\bar{\,\,\, }$ defines a map from the set $\La_m^+(r)$ of all $m$-partitions
of $r$ to $\La(r)$. 

For each $m$-composition $\bla =(\la^{(1)}, \cdots, \la^{(m)})$,   
define $[\bla]=[a_0, a_1, \cdots, a_m]$ such that $a_0=0$ and $a_i=
\sum_{j=1}^i |\la^{(j)}|$. Recall from \cite{DR, \S1} that
the set of all $[\bla]$ form a poset
$\La[m,r]$ which is isomorphic to the poset $\La(m,r)$ of all
compositions of $r$ with $m$ parts. Here the partial ordering
on $\La[m,r]$ is given by $\lecq$:
 $[a_i]\lecq [b_i]$ if $a_i\le b_i$, $i=1, \cdots, m$, while
$\La(m,r)$ has the usual dominance order $\trianglelefteq$.

For any $\bka=[a_i]\in\La[m,r]$, following \cite{DR, (1.6)}, let 
 $w_{\boldkey a}\in \frak S_r$  be the element defined by 
$$
(a_{i-1}+l) w_{\boldkey a}=r-a_i+l\text{ for all }i
\text{ with } a_{i-1}<a_i, 1\le l\le a_{i}-a_{i-1}.
$$
For example,
$w_{\bka}=\left(\smallmatrix 1\,\,2\,\,3\,\,4 &5\,\,6\,\,7\,\,8 & 9\cr
                  \underline{6\,\,7\,\,8\,\,9}&\underline{2\,\,3\,\,4\,\,5}
&\underline{1}\cr
\endsmallmatrix\right),$
if $\bka=[0,4,8,9]$.
Note that $w_\bka$ is the unique distinguished double coset representative
satisfying 
$$w_\bka^{-1}\fS_\bka w_\bka=\fS_{\bka'},\tag1.3$$
where, for $\bka=[a_0,a_1,\cdots,a_m]\in\La[m,r]$, 
$$\bka'=[r-a_m,r-a_{m-1},\cdots, r-a_0]$$
(see \cite{DR, (1.2)}),
and $\fS_\bka$ and $\fS_{\bka'}$ are the Young subgroups defined in
 (1.1).

We may also identify
$\bla$ with its Young diagram $\sY(\bla)$. For example,
$\bla=((31), (22), (1))$ is identified with
$$\sY(\bla)=\boxed{
\smallmatrix 
\square&\square&\square&\hskip.3cm\square&\square&\hskip.3cm\square\cr
\square&       &       &\hskip.3cm\square&\square&\hskip.3cm       \cr
\endsmallmatrix}.$$
Let $\bft^\bla$ be the $\bla$-tableau obtained
by putting the numbers $1,\cdots,r$ in order into
the  boxes in the Young diagram of $\bla$
from left to right down successive rows.
From the example above, we have
$$
\bft^\bla
=\boxed{\smallmatrix 1 & 2& 3& \hskip.3cm 5 & 6  &\hskip.3cm 9\cr
                 4 &  &  & \hskip.3cm 7 & 8  &\hskip.3cm   \cr
                  \endsmallmatrix}.$$
We also define the $\bla$-tableau $\bft_\bla$ by putting the numbers
from right to left down successive columns as illustrated in the following 
example
$$
\bft_{\bla}=\boxed{\smallmatrix 6 & 8& 9& \hskip.3cm 2 & 4  &\hskip.3cm 1\cr
                 7 &  &  & \hskip.3cm 3 & 5  &\hskip.3cm   \cr
                  \endsmallmatrix}.$$
Now, associated to an $m$-partition $\bla =(\la^{(1)}, \cdots, \la^{(m)})$
of $r$, we define the element $w_\bla\in\fS_r$ by $\bft^\bla w_\bla
=\bft_\bla$. More precisely,  if $\bft^i$  (resp. $\bft_i$) denote the
$i$-th subtableau of $\bft^\bla$ (resp. $\bft_\bla  w_{[\bla]}^{-1}$) and define
$w_{(i)}$ by $\bft^iw_{(i)}=\bft_i$, then  
$\bft^\bla w_{{(1)}}\cdots w_{{(m)}}w_{[\bla]}=\bft_\bla $. Likewise,
if we define $\tilde \bft^i$ (resp. $\tilde \bft_i $) the $i$-th
subtableau of $\bft^\bla  w_{[\bla]}$ (resp. $\bft_\bla$) and  $\tilde
w_{(i)}$ with  $\tilde \bft^i \tilde w_{(i)}=\tilde
 \bft_i$, then   
$\bft^\bla w_{[\bla]}\tilde w_{(1)} \cdots \tilde w_{(m)}=\bft_\bla$. 
We have, therefore
$$w_\bla
= w_{{(1)}}\cdots  w_{(m)}  w_{[\bla]}
=w_{[\bla]}\tilde w_{{(m)}}\cdots\tilde w_{{(1)}},\quad
w_{[\bla]}^{-1} w_{{(i)}}w_{[\bla]}=\tilde w_{(m-i+1)}.\tag1.4
$$
Note that $w_{(i)}w_{(j)}=w_{(j)}w_{(i)}$ and  $\tilde w_{(i)}\tilde  
w_{(j)}= \tilde w_{(j)} \tilde w_{(i)}$ for $ i, j=1, 2, \cdots, m$.
Clearly, by (1.2-3), we have the following trivial intersection property
$$w_\bla^{-1}\fS_{\bar\bla} w_\bla\cap \fS_{\bar{\bla'}}=\{1\}.\tag1.5$$

The usual dominance order on $\La^+(r)$ can also be generalized to
$\La_m^+(r)$.  For $m$-partitions
$\bla=(\la^{(1)}, \cdots, \la^{(m)})$ and 
$\bmu=(\mu^{(1)}, \cdots, \mu^{(m)})$ with
$[\bla]=[a_0,\cdots,a_m]$ and $[\bmu]=[b_0,\cdots,b_m]$, we say that
$\bla\trianglelefteq \bmu$
if 
$$
\sum_{j=1}^i a_j +\sum_{k=1}^l \la^{(i+1)}_k
\le \sum_{j=1}^i b_j +\sum_{k=1}^l \mu^{(i+1)}_k
$$
for all $i=0,1, \cdots, m-1$ and all $k$.
The following lemma follows immediately from the definition.

\proclaim{(1.6) Lemma} Suppose $\bla, \bmu\in \La_m^+(r)$ and
 $[\bla]=[\bmu]$. Then 
$\bla\trianglelefteq \bmu$ if and only if $\la^{(i)}\trianglelefteq 
\mu^{(i)}$ for all $i$.
\endproclaim

Let $W=(\BZ/m\BZ)\wr\fS_r$ be
the wreath product of the cyclic group of order $m$ and
 the symmetric group $\fS_r$. Then $W$ is not a Coxeter group
except for $m=1,2$. 
 The Ariki-Koike algebra $\bH_R=\bH^r_{mR}$ 
over a {\it commutative ring $R$ with 1} is an
Iwahori-Hecke type algebra associated to $W$ and
parameters $q,q^{-1},u_1,\cdots,u_m\in R$: it is
associative with generators $T_0,T_1,\cdots,T_{r-1}$ 
subject to the  relations:  
$$\cases (T_0-u_1)\cdots (T_0-u_m)=0,&\cr
         T_0T_1T_0T_1 = T_1T_0T_1T_0,& \cr 
         (T_i-q)(T_i+1)=0, &1\le i\le r-1\cr 
         T_iT_{i+1}T_{i}=T_{i+1}T_iT_{i+1}, &1\le i\le r-2\cr 
         T_i T_j=T_jT_i,&0\le i<j-1\le r-2.\cr 
\endcases \tag 1.7$$
Note that $\bH_R$ is the group algebra of $W$ when $q=1$ and $u_1=\xi$ is a
primitive $m$-th root of unity and $u_i=\xi^i$.  
Let $\sH_R=\sH_R(r)$ denote the subalgebra of $\bH_R$ generated by $T_1,\cdots,T_{r-1}$.
Then $\sH_R=\bH^r_{1R}$ is the Iwahori-Hecke algebra associated to $\fS_r$.
We will denote by $\sH_R(\fS_\la)$ 
the subalgebra associated to the Young subgroup $\fS_\la$ of $\fS_r$.
Let
$$x_\la=\sum_{w\in\fS_\la}T_w \text{ and }y_\la =
\sum_{w\in\fS_\la}(-q)^{-l(w)}T_w.\tag1.8$$
Then $x_\la\sH_R$ and $y_\la\sH_R$ define the associated $q$-permutation modules,
and the Specht module \cite{DJ} associated to a partition $\la$ of $r$
 is the {\it submodule}
$S^\la=x_\la T_{w_\la}y_{\la'}\sH_R$ of $x_\la\sH_R$.

We now come to define the notion of Specht modules for $\bH^r_{mR}$.
Let $L_1=T_0$, $L_i=q^{-1} T_{i-1} L_{i-1} T_{i-1} $ for 
$i=2, \cdots, r$, and put $\pi_0 (x)=1$, 
$\pi_{a} (x)= \prod_{j=1}^{a} (L_j -x)$ for any  $x\in R$ and any positive 
integer $a$. We define, for $\bka=[\bla]=[a_0,a_1,\cdots,a_m]\in\La[m,r]$, 
$$
\pi_{\bka}=\pi_{a_1}(u_2)   \cdots \pi_{a_{m-1}}(u_m)
\text{ and }
\tpi_{\bka} =\pi_{a_1} (u_{m-1})   \cdots \pi_{a_{m-1}} (u_1),
$$
and, for  $\bla\in\La_m^+(r)$, 
 $$x_\bla=\pi_{[\bla]} x_{\bar\bla}=x_{\bar\bla}\pi_{[\bla]}\quad\text{ and }\quad
y_\bla=\tpi_{[\bla]} y_{\bar\bla}=y_{\bar\bla}\tpi_{[\bla]} .$$

\proclaim{(1.9) Lemma} For any $m$-partition  
 $\bla=(\la^{(1)}, \cdots, \la^{(m)})$,  the $R$-submodule
 $x_\bla \bH y_{\bla'}=Rz_\bla$ is free of rank 1, where
$z_\bla=x_\bla
 T_{w_\bla} y_{\bla'}$. 
\endproclaim

\demo{Proof} From \cite{DR, (2.8), (3.1)},
 we have $x_\bla \bH y_{\bla'}
=x_{\bar \bla}   \sH(\frak S_{\bka})v_\bka y_{\bar{\bla'}} $ where
 $\bka=[\bla]$ and $v_\bka =\pi_\bka T_{w_\bka}
\tpi_{\bka'}$.  Since 
$y_{\bar{\bla'}}=y_{\la^{(m)'}\vee \cdots \vee \la^{(1)'}}$
and $ v_\bka y_{\bar{\bla'}}= 
y_{\la^{(1)'}\vee \cdots \vee \la^{(m)'}}v_\bka$, 
we have,  by \cite{DJ, (4.1)},  the $R$-submodule 
$ x_{\la^{(1)}\vee \cdots\vee \la^{(m)}}\sH(\frak
 S_{\bka})y_{\la^{(1)'}\vee \cdots \vee \la^{(m)'}}$ is generated by
the element 
$
Z_{\bar\bla }= x_{\la^{(1)}\vee \cdots \vee \la^{(m)} }T_{ w_{(1)}
\cdots w_{(m)}}
y_{\la^{(1)'}  \vee \cdots \vee \la^{(m)'}}$.  
 Therefore, by \cite{DR, (3.4)}, the 
$R$-submodule  $x_\bla \bH y_{\bla'}$ is 
of rank  1, generated by $z_\bla=x_\bla
 T_{w_\bla} y_{\bla'}=Z_{\bar\bla }v_\bka$.\qed 
\enddemo

\definition{ (1.10) Definition} The right $\bH_R$-module
$$S^\bla= x_\bla \bH_R y_{\bla'} \bH_R=
 x_\bla T_{w_\bla}y_{\bla'}\bH_R
$$
is called the {\it Specht module} associated to the $m$-partition $\bla$. 
Similarly, the right $\bH_R$-module  
$$
\tilde S^\bla=y_\bla \bH_R  x_{\bla'}\bH_R =
y_\bla T_{w_\bla}x_{\bla'}\bH_R$$
is called the {\it twisted Specht module} of $\bH_R$ associated to the
$m$-partition
$\bla$.
\enddefinition
 The name will be justified in (2.5)
since $\{S^\bla\mid \bla\in \La_m^+(r)\}$ is 
a complete set of non-isomorphic simple modules for a semi-simple
$\bH_R$ \cite{DR, (5.1)}. 

\subhead 2. The Standard Basis Theorem \endsubhead
The standard basis theorem for Hecke algebra of type $A$
refers to a basis for a Specht module $S^\la$ indexed by the
standard $\la$-tableaux. Recall from \cite{DJ, \S1}, a
 $\la$-tableau of the from
$\bft^\la w$, $w\in\fS_r$, is called  standard 
if it has increasing rows and columns.
The following theorem is due to Dipper and James \cite{DJ, (5.6)}.

\proclaim{(2.1) Theorem}
Let $\la$ be a partition of $r$.
 Then the Specht module $S^\la=z_\la\sH_R$ with
$z_\la=x_\la T_{w_\la}y_{\la'}$ is a free $R$-submodule 
with basis
$$\{z_\la T_d\mid d\in\fS_r,\bft_\la d \text{ is standard}\}
=\{z_\la T_d\mid d\in\fS_r,\bft^{\la'} d \text{ is standard}\}.$$ 
\endproclaim

We are now ready to generalize this result to Ariki-Koike algebras.
Let $\bla$ be an $m$-partition of $r$. A $\bla$-tableau $\bft$ 
is said to be  {\it standard} if each $\la^{(i)}$-subtableau 
has increasing rows and columns.  

\proclaim{(2.2) Standard Basis Theorem} Let $\bla$ be an 
$m$-partition of $r$. Then  the Specht module 
$S^\bla=z_\bla\bH_R$ with $z_\bla=x_\bla T_{w_\bla}y_{\bla'}$
 is a free $R$-submodule with basis
$\{z_\bla T_d\mid d\in \fS_r, \bft_\bla d \text{ is standard }\}
=\{z_\bla T_d\mid d\in\fS_r, \bft^{\bla'} d \text{ is standard }
\}$.
\endproclaim

\demo{Proof} We first notice from the proof of (1.9) that
$$
z_\bla = Z_{\bar \bla}v_{[\bla]}\text{ where }
Z_{\bar\bla }= x_{\la^{(1)}\vee \cdots \vee \la^{(m)} }T_{ w_{(1)}
\cdots w_{(m)}}
y_{\la^{(1)'}  \vee \cdots \vee \la^{(m)'}}\tag 2.3
$$ 
Thus, \cite{DR, (3.1)} implies $z_\bla\bH_R=z_\bla\sH_R$.
Let $\la=(|\la^{(1)}|,\cdots,|\la^{(m)}|)$. 
By (2.1), the $\sH_R({\fS_\la})$-module $Z_{\bar\bla }\sH_R({\fS_\la})$
has  basis
$\{Z_{\bar\bla } T_d\mid d\in \fS_\la,
\bft_\bla d \text{  standard }\}$.
Since, as $\sH_R$-module, 
$z_\bla\sH_R\cong Z_{\bar\bla }v_{[\bla]}\sH_R({\fS_\la})
\otimes_{\sH_R({\fS_\la})}\sH_R$ (see \cite{DR, (3.4)}),
 $S^\bla$ has a basis
$\{z_\bla T_dT_w\mid d\in \fS_\la,
\bft_\bla d \text{ standard },w\in\sD_\la\}$.
Note that  $w\in \sD_\la$ if and only if 
 $(a_i+1)w, \cdots,  (a_{i+1})w$ is increasing  for 
$i=0, \cdots, m-1$, where $[a_0,\cdots,a_m]=[\bla]$. Thus, for $d\in \fS_\la$,
$\bft_\bla dw$ is  standard for all $w\in\sD_\la$
if  $\bft_\bla d$  is standard,  and vice-versa. 
Therefore, this  basis is the required basis.
 \qed
\enddemo

\proclaim{(2.4) Corollary}  The module $S^\bla$
(resp. $\tilde S^\bla$) is an $R$-pure submodule of $x_\bla\bH_R$
(resp. $y_\bla\bH_R$).
\endproclaim

\demo{Proof} This follows immediately from the purity of the
type $A$ Specht modules.\qed
\enddemo

\remark{\bf (2.5) Remark} 
The standard basis given in (2.2) can also be obtained by using  the Murphy 
type (or cellular) basis
for $x_\bla\bH_R$. See  \cite{DS, (5.2.1)} for the case $m=2$ and  
\cite{DJM, (4.14)} in general. 
In these work, the bases are all defined for $x_{\bla'}\bH_R$. However, the 
 ring automorphism 
$$
\Phi:\bH_R\rightarrow \bH_R\tag2.6
$$ 
defined by setting 
$\Phi(q)=q^{-1}$, $\Phi(u_i)=u_{m-i+1}$, for $1\le i\le m$,
  $ \Phi (T_j) =-q^{-1}T_{j} $ 
for $1\le j\le r-1$ and $\Phi(T_0)=T_0$ (cf. (1.7)), will turn such a basis 
for $x_{\bla'}\bH_R$ into a basis
for $y_{\bla'}\bH_R$. Denote this basis by 
$$\{Y_{\fku\bft}^\bmu\mid \bmu\trianglerighteq\bla', \bft\in \bT^s(\bmu),\fku
\in\fT^{ss}(\bmu,\bla')\},$$
where $\bT^s(\bmu)$ (resp. $\fT^{ss}(\bmu,\bla')$) is the set of all standard
$\bmu$-tableaux (semi-standard $\bmu$-tableaux of type $\bla'$). 
Note that, since $\Phi(\pi_\bka)=\tpi_\bka$ and 
$\Phi(x_{\bar\bla})=y_{\bar\bla}$, $Y_{\fku\bft}^\bmu$ is of the form $hy_\bmu T_d$ for some
$h\in\sH_R$ and $d\in \fS_r$.

We now consider the homomorphism
$$
\ph:y_{\bla'}\bH_R\to z_\bla\bH_R;\quad y_{\bla'}h\mapsto z_{\bla}h
=x_\bla T_{w_{[\bla]}}(y_{\bla'}h).\tag2.7
$$
we claim that if $\ph(Y_{\fku\bft}^\bmu)=x_\bla T_{w_{[\bla]}}Y_{\fku\bft}^\bmu\neq 0$ then $\bmu=\bla'$.
Using \cite{DR, (2.8)},  
we see easily that $\ph(Y_{\fku\bft}^\bmu)\neq 0$
implies  $[\bla]\lecq [\bmu']$. Since $\bmu\trianglerighteq \bla'$,
we have $[\bmu]\recq [\bla']$ and hence $[\bmu']\lecq [\bla]$. So 
$[\bmu']=[\bla]$. Thus, 
if we write $Y_{\fku\bft}^\bmu=hy_\bmu T_d$ as above, then,
by \cite{DR, (3.1a,c)} and noting $x_{\bar\bla}T_{w_{[\bla]}}=T_{w_{[\bla]}}
x_{\la^{(m)}\vee\cdots\vee\la^{(1)}}$, we have
$x_\bla T_{w_{[\bla]}}Y_{\fku\bft}^\bmu=
v_{[\bla]}x_{\la^{(m)}\vee\cdots\vee\la^{(1)}}h'y_{\bar\bmu}T_d$ for some
$h'\in\sH_R(\fS_{[\bla']})$. Since $[\bmu]=[\bla']$, 
this is non-zero implies, by \cite{DJ, (4.1)}, 
${\la^{(m-i+1)}}\trianglelefteq \mu^{(i)\prime}$ for all $i=1, \cdots, m$. 
On the other hand, by (1.6), $\bmu\trianglerighteq \bla'$ and $[\bmu']=[\bla]$
implies ${\la^{(m-i+1)}}'\trianglelefteq \mu^{(i)}$, $1\le i\le m$. Hence 
${\la^{(m-i+1)}}'= \mu^{(i)}$ for all $i$, and
therefore, 
$\bmu=\bla'$.  This proves that the image of the cellular basis for the 
Specht module defined as the quotient of $y_{\bla'}\bH_R$ modulo the submodule
spanned by all $Y_{\fku\bft}^\bmu$ with $\bmu\triangleright\bla'$
is exactly the standard basis described in (2.2).
 Using a similar argument, we see that the twisted Specht modules
${\tilde S}^\bla$ defined in (1.10) are isomorphic to  those defined
in \cite{DJM, (3.28)} as quotient modules of $x_{\bla'}\bH_R$. 
\endremark

\subhead 3. Branching rules, I\endsubhead
For a partition $\la=(\la_1,\cdots,\la_m)$ of $r$, we identify the
boxes in the Young diagram $\sY(\la)$ with its position coordinates.
Thus, we have
$$\sY(\la)=\{(i,j)\in\BZ^+\times\BZ^+\mid j\le\la_i\}.\tag3.1$$
The elements of $\sY(\la)$ will be called {\it nodes}.
A node of the form $(i,\la_i)$ (resp. $(i, \la_i+1)$) is called {\it
removable} (resp.  {\it addable}) if $\la_i>\la_{i+1}$
(resp. $\la_{i-1}>\la_{i}$). Let
$\bla=(\la^{(1)},\cdots,\la^{(m)})$ be an $m$-partition. Then
 its Young diagram $\sY(\bla)$ is a
union of the Young diagram $\sY(\la^{(k)})$, $1\le k\le m$. Thus, we have
as a set of nodes
$$\sY(\bla)=\{(i,j)_k\mid i,j\in\BZ^+, j\le\la_i^{(k)}, 1\le k\le m\}.$$
A node of $\sY(\bla)$ is said to be removable (resp. addable)
 if it is a removable (resp. addable)
node of $\sY(\la^{(k)})$ for some $k$. Let $\sR_\bla$ denote the set
of all removable nodes of $\sY(\bla)$. Then 
$N=\#\sR_\bla=\sum_{i=1}^m\#\sR_{\la^{(i)}}$, and, if
$\bla'=(\la^{(m)\prime},\cdots,\la^{(1)\prime})$ is the $m$-partition
dual to $\bla$, then we have a bijection
$
\tau:\sR_\bla\to\sR_{\bla'}; (i,j)_k\mapsto (j,i)_{m-k+1}.
$

The ordering on $\sR_\bla$ will be fixed from top to bottom and from left 
to right: $\sR_\bla=\{\rn_1, \cdots,\rn_N\}$,
let $j_\rn$, $\rn\in\sR_\bla$, be the number at the node $\rn$ in 
$\bft_\bla$. Note that, if $[\bla]=[a_i]$ with 
$a_i=\sum_{j\le i}|\la^{(j)}|$, then
$$j_{\rn_{N_i+1}}=r-a_i,\text{ where } N_0=0,\,\,N_i=\#\cup_{j\le i}
\sR_{\la^{(j)}}, 1\le i\le m-1.\tag3.2
$$
For example, for $\bla=((31),(22),(1))$, 
$\sR_\bla=\{(1,3)_1,(2,1)_1,(2,2)_2,(1,1)_3\}$, and the corresponding
$j_\rn$'s are 9,7,5,1.
We have the following result.
Recall for any $i\le j$ the elements
 $s_{i,j}=\left(\smallmatrix i&i+1&i+2&\cdots&j\cr
j&i&i+1&\cdots&j-1\cr\endsmallmatrix\right)$ if $i<j$ and $s_{i,i}=1$.

\proclaim{(3.3) Lemma} Let $\rn\in\sR_\bla$, and let $\bla_\rn$ denote the 
multi-partition associated to the tableau obtained by 
removing $\rn$ from $\sY(\bla)$. Let $\bT^s(\bla)$ be the set
of all standard $\bla$-tableaux. Then:

(a) The $\bla_\rn$-tableau $\bft_{\bla_\rn}$ is the same as the tableau
obtained by removing the entry $r$ from $\bft_\bla s_{j_\rn,r}$.
Likewise, the $\bla_\rn'$-tableau $\bft^{\bla_\rn'}$ is the tableau obtained 
by removing the entry $r$ from $\bft^{\bla'} s_{j_\rn,r}$.

(b) $w_\bla =s_{i_\rn,r} w_{\bla_\rn}s_{r,j_\rn}$,  where $i_\rn$ is defined by
$(i_\rn)w_\bla=j_\rn$.

(c) $s_{j_\rn,r}\in\sD_{\bar{\bla'},(r-1,1)}$
and $\fS_{\bar{\bla'}}^{s_{j_\rn,r}}\cap\fS_{r-1}=\fS_{\bar{\bla_\rn'}}$.

Moreover, we have 
$$\{d\in\fS_r\mid \bft_\bla d\in\bT^s(\bla)\}=\bigcup_{\rn\in\sR_\bla}
\{s_{j_{\rn},r}x\mid x\in\fS_{r-1},\bft_{\bla_{\rn}}x\in\bT^s(\bla_{\rn})\}.
$$
\endproclaim

\demo{Proof}  The statement (a) follows
immediately from the definition of $s_{i,r}$.
Since $j_\rn$ is the number at $\rn$ in $\bft_\bla$ and
$i_\rn$ is the number at $\rn$ in $\bft^\bla$
by (a), we see that $\bft^{\bla_\rn}$ can be obtained by removing the entry 
$r$ at $\rn$
from $\bft^\bla s_{i_\rn,r}$. 
On the other hand, since $\bft^\bla w_\bla=\bft_\bla$, 
it follows that
$(\bft^\bla s_{i_\rn,r})s_{r,i_\rn}w_\bla s_{j_\rn,r}=
\bft^\bla w_\bla s_{j_\rn,r}=\bft_\bla s_{j_\rn,r}$. Therefore, $s_{r,i_\rn}
w_\bla s_{j_\rn,r}=w_{\bla_\rn}$, proving (b).
The first assertion in (c) follows from the relation 
$$
s_{r, j_\rn} s_i s_{ j_\rn, r}=\cases s_i,  & \text{ if $i<j_\rn$,}\cr
                                      s_{i-1},  & \text{ if $i>j_\rn$,}\cr
        s_{r, j_\rn+1} s_{j_\rn}s_{ j_\rn+1, r} & \text{ if $i=j_\rn$.}\cr 
\endcases
$$
For $d$ with $\bft_\bla d\in\bT^s(\bla)$, write $d=d_1x$ with
$x\in\fS_{r-1}$ and $d_1$ distinguished relative to $\fS_{r-1}$. Then,
$d_1=s_{i,r}$ for some $i$. Since $\bft_\bla d_1$ is again standard,
it forces $i=j_\rn$ for some $\rn\in\sR_\bla$. 
Thus,  $\bft_{\bla_\rn}$ is obtained by removing the entry $r$
from $\bft_\bla d_1$. Therefore,
the tableau $\bft_{\bla_\rn}x$ can be obtained from  
$\bft_{\bla}s_{j_\rn, r} x$ by removing the entry $r$. Since 
$\bft_{\bla}s_{j_\rn, r} x$ is standard,  $\bft_{\bla_\rn}x$ is
standard, too, proving the inclusion ``$\subseteq$''.
Since $\bft_\bla s_{j_\rn, r}x$
is the tableau obtained from $\bft_{\bla_\rn}x$ by adding the entry
$r$ at the node $\rn$,   $\bft_\bla s_{j_\rn, r}x$ is a
standard tableau, proving the inclusion  ``$\supseteq$''. Therefore,
the required equality holds. \qed\enddemo

\proclaim{(3.4) Lemma} (a) For positive integers $j,k$ with
$1\le j, k\le r-1$, we have  
$$
T_{j,r}L_k=\cases L_kT_{j,r},&\text{ \rm if }k<j\cr
L_{k+1}T_{k,r}-(q-1) L_{k+1}T_{k+1,r},
&\text{ \rm if }k=j\cr
L_{k+1}T_{j,r}-(q-1) L_{k+1}T_{k+1, r}T_{j,k},&\text{ \rm if }k>j,
\endcases
$$
where $T_{i,l}=T_{s_{i,l}}$.

For $\bka=[0, a_1, \cdots, a_{m-1}, r]\in\La[m, r]$, let $v_\bka=\pi_\bka
T_{w_{\bka}}\tpi_{\bka'}$.   Then  

(b) $v_\bka T_i=T_{(i)w_\bka^{-1}} v_\bka$ if $i\neq r-a_j$ for $j=1,
\cdots, m$;

(c) $v_{\bka} L_k= u_{j} q^{r-a_{j}+1-k}v_\bka T_{k, r-a_j+1} T_{r-a_j+1, k}$,
where $j$ satisfies $r-a_{j}+1\le k\le r-a_{j-1}$.
\endproclaim

\demo{Proof} For $k<j$, the result in (a) is clear. If $k=j$, then 
$T_{k,r}L_k=T_kL_kT_{k+1,r}=q L_{k+1}T_k^{-1}T_{k+1,r};$
if $k>j$, then
$$\aligned
T_{j,r}L_k
&=T_{j,k+1}L_kT_{k+1,r}=T_{j,k}T_kL_kT_{k+1,r}\cr
&=qT_{j,k} L_{k+1}T_{k}^{-1} T_{k+1,r}=q  L_{k+1}T_{j,k}
  T_k^{-1} T_{k+1,r}.\cr
\endaligned$$
Now the result (a) follows easily since $qT_k^{-1}=T_k-(q-1)$.

The statements (b) and (c) have been proved in \cite{DR, (3.1)}.\qed
\enddemo

The Branching Theorem for symmetric groups can be found in \cite{JK}.
The following is the $q$-version (see \cite{Jo, 3.4}).

\proclaim{(3.5) Lemma}  
For any $\la\in\La^+(r)$, let
$\rn_1,\cdots,\rn_k$
be the removable nodes of $\sY(\la)$ 
counted from top to bottom, and  define $M_0=0$ and
$M_t=z_\la T_{j_{\rn_t},r}\sH_R(\fS_{r-1})+M_{t-1}$ for $t\ge 1$.
Then we have a filtration of $\sH_R(\fS_{r-1})$-submodules for $S_R^\la=z_\la\sH_R$:
$$0=M_0\subset M_1\subset\cdots\subset M_k=S_R^\la$$
with sections of Specht modules:
$M_t/M_{t-1}\cong S_R^{\la_{\rn_t}}$.
\endproclaim

We are now ready to prove the (ordinary) Branching Theorem for Ariki-Koike
algebra. A version of this result has been given by Ariki-Mathas
\cite{AM, 1.5} with a proof using cellular basis.
If we turn their result into a result for $y_{\bla'}\bH$ using (2.6) and
apply the homomorphism $\ph$ given in (2.7), the theorem below
follows immediately from theirs. 
However, the proof supplied here can be viewed as a proof in our context
without using cellular bases, following the idea of lifting from
level $1$ to level $m$.
Recall from (2.2) that $z_\bla=x_\bla T_{w_\bla} y_{\bla'}$.

\proclaim{(3.6) Branching Theorem}  
Let $\bla$ be an $m$-partition of $r$, and let 
$\rn_1,\cdots,\rn_N$ be the removable nodes of $\sY(\bla)$
counted from top to bottom and from left to right.
Define $\bkM_0=0$ and
$\bkM_t=z_\bla T_{j_{\rn_t},r}\bH^{r-1}_{mR}+\bkM_{t-1}$ for $t\ge 1$.
Then these modules form a filtration of $\bH^{r-1}_{mR}$-submodules for $
S_R^\bla=z_\bla\bH_R$:
$$0=\bkM_0\subset \bkM_1\subset\cdots\subset \bkM_N=S_R^\bla$$
with sections of Specht modules:
$\bkM_t/\bkM_{t-1}\cong S_R^{\bla_{\rn_t}}$.
\endproclaim
\demo{Proof}
Let 
$M_0=0$ and
$M_t=z_\bla T_{j_{\rn_t},r}\sH_R(\fS_{r-1})+M_{t-1}$ for $t\ge 1$.
We claim that $\bkM_t=M_t$ for all $t$.  
Indeed, since  $M_t$ is a right $\sH_R(\fS_{r-1})$-module, we need
only prove 
$$
z_\bla T_{j_{\rn_t},r}L_k\in M_t,\quad\forall k\le r-1.\tag3.7
$$
Let $[\bla]=[a_0,a_1,\cdots,a_m]$, and choose $i$ such that 
$r-a_{i+1}< j_{\rn_t}\le  r-a_{i}$. Thus, by (3.2),  $\rn_t$ is a
removable node of the 
$(i+1)$-th Young subdiagram $\sY(\la^{(i+1)})$.
By (2.3) and (3.4b), we have
$$
z_\bla T_{j_{\rn_t},r}L_k =Z_{\bar\la}v_{[\bla]}T_{j_{\rn_t},r}L_k=z_{\bar\la}  
T_{(j_{\rn_t})w_{[\bla]}^{-1}}\cdots 
T_{(r-a_{i}-1)w_{[\bla]}^{-1}}v_{[\bla]} T_{r-a_i,r}L_k.
$$
Put $b_j=r-a_j$ for all $j$. 
By (3.4a,c), we have the following equalities:
$$
z_\bla T_{j_{\rn_t},r}L_k
=\cases u_j q^{b_j+1-k} z_\bla
T_{j_{\rn_t}, r} T_{k, b_j+1} T_{b_j+1, k}, &
\text{ if }k<b_i,b_j<k\le b_{j-1}\cr
u_i z_\bla T_{j_{\rn_t}, r}-h_1, &
\text{ if }k=b_i\cr
u_jq^{b_j-k} z_\bla  T_{j_{\rn_t}, r} T_{k, b_j}T_{b_{j}, k}-h_2,
& \text{ if }k>b_i,b_j\le k< b_{j-1}
\endcases \tag 3.8
$$
where 
$$\aligned 
           & h_1=u_i(q-1) z_{\bla} T_{j_{\rn_t}, b_i} T_{b_i+1,r}
                = u_i(q-1) z_{\bla}  T_{b_i+1,r} T_{j_{\rn_t}, b_i},\cr 
           & h_2=u_jq^{b_j-k}(q-1) z_\bla T_{k+1, b_j+1} T_{b_j+1, r} 
                 T_{j_\rn, k} .\cr\endaligned 
$$
Since 
$ T_{k+1, b_j+1} T_{b_j+1, r}\in \sum_{l\ge k>j_\rn} T_{l, r}\sH_R(\fS_{r-1})$,
and $j_{\rn_1}, \cdots, j_{\rn_N}$ is  decreasing, 
 we have $h_1, h_2\in M_{t-1}$ by induction, proving  (3.7), 
and hence the claim.     

Put $z_\bla=v_{[\bla]}\tilde Z_{\bar{\bla^\circ}}$ (cf. (2.3)), 
where $\bla^\circ=(\la^{(m)},\cdots,\la^{(1)})$ and
$$\tilde Z_{\bar{\bla^\circ}}=
\tilde Z_{\la^{(m)}\vee\cdots\vee\la^{(1)}}=x_{\la^{(m)}\vee\cdots\vee\la^{(1)}}T_{\tilde w_{(m)}
\cdots\tilde w_{(1)}}y_{\la^{(m)\prime}\vee\cdots\vee\la^{(1)\prime}}.$$
Note that $y_{\la^{(m)\prime}\vee\cdots\vee\la^{(1)\prime}}=y_{\bar{\bla'}}$.
Since $\tilde w_{(i)}\in
\fS_{\{r-a_i+1,\cdots,r-a_{i-1}\}}$ (see the definition above (1.4))
and, for $a<b$, $\fS_{\{a+1,\cdots,b\}}^{s_{a,b}}\cap \fS_{\{a,\cdots,b-1\}}=
 \fS_{\{a,\cdots,b-1\}}$,
(3.3c) for $m=1$ implies that $y_{\bar{\bla'}}T_{j_{\rn_t},r}=
h_0T_{j_{\rn_t},r}y_{\bar \bla_{\rn_t}'}$
for some $h_0\in\sH_R(\fS_{\{r-a_{i+1}+1,\cdots,r-a_i\}})$. 
On the other hand,
 by (3.3b), 
$$\aligned
&\quad\,T_{\tilde w_{(m)}\cdots\tilde w_{(1)}}h_0T_{j_{\rn_t},r}\cr
&=T_{\tilde w_{(m)}\cdots\tilde w_{(i+2)}}
(T_{\tilde w_{(i+1)}}h_0T_{j_{\rn_t},r-a_i})
(T_{\tilde w_{(i)}}T_{r-a_i,r-a_{i-1}})\cdots
(T_{\tilde w_{(1)}}T_{r-a_1,r})\cr
&=T_{\tilde y_{(m)}\cdots\tilde y_{(i+2)}}
[(T_{\tilde w_{(i+1)}}h_0T_{j_{\rn_t},r-a_i})
T_{r-a_i,r}]T_{\tilde y_{(i)}\cdots\tilde y_{(1)}},\cr
\endaligned
$$
 where $\tilde y_{(i)}$ is the $\tilde w_{(i)}$ defined
 relative to $\bla_{\rn_t}$ as in (1.4),
and  
$x_{\la^{(m)}\vee\cdots\vee\la^{(1)}}T_{r-a_i,r}
=x_{\la^{(m)}\vee\cdots\vee\la^{(i+1)}}T_{r-a_i,r}
\tilde x_{\la^{(i)}_{\rn_t}\vee\cdots\vee\la^{(1)}_{\rn_t}}$, where 
$\tilde x_{\la^{(i)}_{\rn_t}\vee\cdots\vee\la^{(1)}_{\rn_t}}$ is the
sum of $T_w$'s with $w\in\fS_{\{r-a_{i},\cdots,r-a_{i-1}-1\}}\times\cdots\times
\fS_{\{r-a_1,\cdots,r-1\}}$.
So, noting $\la^{(j)}=\la_{\rn_t}^{(j)}$ for all $j\neq i+1$, we eventually see that 
$$
Z_{\bar{\bla^\circ}}T_{j_{\rn_t},r}
=\tilde Z_{\la^{(m)}_{\rn_t}\vee\cdots\vee\la^{(i+2)}_{\rn_t}}
[(\tilde x_{\la^{(i+1)}}T_{\tilde w_{(i+1)}}\tilde y_{\la^{(i+1)\prime}}T_{j_{\rn_t},r-a_i})
T_{r-a_i,r}]
{\tilde Z}_{\la^{(i)}_{\rn_t}\vee\cdots\vee\la^{(1)}_{\rn_t}},
$$
where $\tilde x$ and $\tilde y$ are defined over 
$\fS_{\{r-a_{i+1}+1,\cdots,r-a_i\}}$ and the last $\tilde Z$
is defined over $\fS_{\{r-a_{i},\cdots,r-a_{i-1}-1\}}\times\cdots\times
\fS_{\{r-a_1,\cdots,r-1\}}$. 
By Lemma (3.5), we have an epimorphism from
$[(\tilde x_{\la^{(i+1)}}T_{\tilde w_{(i+1)}}\tilde y_{\la^{(i+1)\prime}}T_{j_{\rn_t},r-a_i})
T_{r-a_i,r}]
\sH_R(\fS_{\{r-a_{i+1}+1,\cdots,r-a_i-1\}})$ onto
$\tilde z_{\la^{(i+1)}_{\rn_t}}\sH_R(\fS_{\{r-a_{i+1}+1,\cdots,r-a_i-1\}})$,
where 
$\tilde z_{\la^{(i+1)}_{\rn_t}}=\tilde x_{\la^{(i+1)}_{\rn_t}}
T_{\tilde y_{(i+1)}}
\tilde y_{\la^{(i+1)}_{\rn_t}}$  is defined over $\fS_{\{r-a_{i+1}+1,\cdots,r-a_i
-1\}}$. This results in  an epimorphism  of $\sH_R(\fS_{[\bla_{\rn_t}^\circ]})$-modules
$\tilde Z_{\bar\bla^\circ}T_{j_{\rn_t},r}
\sH_R(\fS_{[\bla_{\rn_t}^\circ]}) \twoheadrightarrow\tilde Z_{\bar\bla_{\rn_t}^\circ}\sH_R(\fS_{[\bla_{\rn_t}^\circ]})
$, and
hence, an  epimorphism  of $\sH_R(\fS_{r-1})$-modules
$\tilde Z_{\bar\bla^\circ}T_{j_{\rn_t},r}
\sH_R(\fS_{{r-1}})  \twoheadrightarrow\tilde Z_{\bar\bla_{\rn_t}^\circ}\sH_R(\fS_{{r-1}})
$.
Observing from (3.5) the kernel of the first epimorphism above, we see
the kernel of this  epimorphism is the intersection
$\tilde Z_{\bar\bla^\circ}T_{j_{\rn_t},r}
\sH_R(\fS_{{r-1}}) \cap N_{t-1}$,
where $N_i$ is defined as $N_i=\tilde Z_{\bar\bla^\circ}T_{j_{\rn_i},r}
\sH_R(\fS_{{r-1}})+N_{i-1}$ for $i\ge 1$, and $N_0=0$.
(Thus, $M_i=v_{[\bla]}N_i$ for all $i$.)
Since the left multiplications by $v_{[\bla]}$ and $v_{[\bla_{\rn_t}]}$
induce isomorphisms by \cite{DR, (3.4)}, it yields an epimorphism
of  $\sH_R(\fS_{r-1})$-modules
$z_\bla T_{j_{\rn_t},r}\sH_R(\fS_{r-1})\twoheadrightarrow S^{\bla_{\rn_t}}$ 
with kernel $z_\bla T_{j_{\rn_t},r}\sH_R(\fS_{r-1})\cap M_{t-1}$.
Therefore,  
we finally obtain an $\sH_R(\fS_{r-1})$-module isomorphism
$f_t:S^{\bla_{\rn_t}}\to M_t/M_{t-1}$ which maps
$z_{\bla_{\rn_t}}h$ to $z_\bla T_{j_{\rn_t},r}h+M_{t-1}$ for
all $h\in\sH_R(\fS_{r-1})$. 

By the claim, it remains to prove that $f_t$ is an
$\bH^{r-1}_{mR}$-module isomorphism.
We first note from (3.4b,c) that (3.8) holds if
 $z_\bla T_{j_{\rn_t},r}$ is replaced by $z_{\bla_{\rn_t}}$ and
$h_i$ by zeros.
Also, note that, for any $x\in\fS_{r-1}$, $T_xT_0=L_kh$ for $h\in\sH_R(\fS_{r-1})$.
So we have for some $h'\in\sH_R(\fS_{r-1})$
$$
f_t(z_{\bla_{\rn_t}}T_xT_0)=f_t(z_{\bla_{\rn_t}}L_k h)
=f_t(z_{\bla_{\rn_t}}h' h)=f_t(z_{\bla_{\rn_t}})h' h
=f_t(z_{\bla_{\rn_t}}T_x)T_0.
$$
Therefore, $f_t$ is an $\bH^{r-1}_{mR}$-module isomorphism, and the theorem is
proven.\qed
\enddemo

\proclaim{(3.9) Corollary} Let $R$ be a field and assume that
$\bH^r_{mR}$ is semi-simple. Then we have an isomorphism of 
$\bH^{r-1}_{mR}$-modules:
$S^\bla_R\vert_{\bH^{r-1}_{mR}}\cong\oplus_{\rn\in \sR_\bla}S^{\bla_\rn}_R$.
\endproclaim


\subhead 4. Branching rules, II\endsubhead
In this section, we shall describe the modular branching theorem
 for an Ariki-Koike algebras $\bH=\bH_{mR}^r$ over a
{\it field} $R$ which has a semi-simple ``bottom'', that is, satisfying
the assumption  
$$f_{m, r}=\prod_{i=1}^{m-1} \prod_{j=i+1}^m  \prod_{k=1-r}^{r-1}
   (u_i q^k-u_j)\neq 0.\tag4.1$$

 Let $l$ be the smallest integer $a$ such that 
$1+q+\cdots+q^{a-1}=0$. (Note that, if $q=1$, then $l$ is the 
characteristic of $R$.)
If such an integer does not exist, then we set 
$l=\infty$.  A partition of $r$ is said to be
$l$-regular if $\la$ has {\it no} non-zero part occurring $l$ or more times.
An $m$-partition $\bla=(\la^{(1)}, \cdots, \la^{(m)})$ is said to be
$l$-regular if each $\la^{(i)}$ is $l$-regular.
Unless otherwise specified, we shall assume in this section
that $R$ is a {\it field} containing the element
$u_1, \cdots, u_m$ and $q\neq 0$, and the subscript $R$ in $\bH_R$ etc. will 
be dropped
for notational simplicity.

\proclaim{(4.2) Theorem} Let $\bH$ be the Ariki-Koike algebra
over a field $R$ such that $f_{m, r}\neq 0$. Then 
the set $\{ y_{\bla'} T_{w_{\bla'}} z_\bla\bH
\mid \bla\in\La_m^+(r) \text{ $l$-regular}\}$ is a complete
set of non-isomorphic simple $\bH$-modules.\endproclaim

\demo{Proof} With our hypothesis, recall from \cite{DR, (3.9),(4.14)} that
 there exist orthogonal idempotents $e_\bka=v_\bka T_{w_{\bka'}}c_\bka$ of $\bH$, (Note that $c_\bka$ is denoted $z_\bka$ in \cite{DR}.)
$\bka\in\La[m,r]$,
such that the categories of $\bH$-modules and
$\epsilon \bH \epsilon$-modules  are Morita equivalent, where 
$\epsilon=\sum_{\bka\in \La[m, r]} e_\bka$.
Let $D^\bla=y_{\bla'} T_{w_{\bla'}} x_{\bla} T_{w_\bla} y_{\bla'}\bH$.
By the Morita equivalence, it suffices to prove that the set
$$
\{D^{\bla}\epsilon\mid \bla\in\La_m^+(r) \text{ $l$-regular}\}
$$
is a complete set of non-isomorphic simple $\epsilon \bH \epsilon$-modules.

Suppose $\bla=(\la^{(1)}, \cdots, \la^{(m)})$ is an
$l$-regular $m$-partition and write $\bka=[a_i]=[\bla]$.
Then, by (2.3), $y_{\bla'} T_{w_{\bla'}} z_\bla
=y_{\bla'} T_{w_{\bla'}}Z_{\bar\bla}v_\bka$. Since 
$v_\bka\bH=e_\bka\bH$, we have by \cite{DR, (3.10),(3.8)}
$$\aligned
D^{\bla}\epsilon
&=y_{\bla'} T_{w_{\bla'}}Z_{\bar\bla}(e_\bka\bH e_\bka)\cr
&=\tpi_{\bka'} T_{{w_\bka}^{-1}}
(y_{{\la^{(1)}}'\vee\cdots\vee {\la^{(m)}}'} T_{
w(1)^{-1}\cdots w(m)^{-1} } Z_{\bar\bla})(e_\bka\bH e_\bka)\cr
&=\tpi_{\bka'} T_{{w_\bka}^{-1}}D^{\bar\bla}e_\bka,\cr
\endaligned\tag4.3$$
where 
$D^{\bar\bla}=(y_{{\la^{(1)}}'\vee\cdots\vee {\la^{(m)}}'} T_{
w(1)^{-1}\cdots w(m)^{-1} } Z_{\bar\bla})\sH(\fS_{\bka})$. 
 Now, it is known  (\cite{J, 8.1i}), the $\sH(\fS_{\{a_{i-1}+1, \cdots,
a_i\}})$-module
$D^{\la^{(i)}}= 
y_{\la^{(i)'}}T_{ w_{(i)}^{-1}} x_{\la^{(i)}} T_{ w_{(i)}}
y_{\la^{(i)'}}\sH(\fS_{\{a_{i-1}+1,\cdots, a_i\}})$ is simple 
 for every $i=1, \cdots, m$, where  $ x_{\la^{(i)}}$
and $y_{\la^{(i)}}$ are defined as in (1.8) using the subgroup
$\fS_{\{a_{i-1}+1,\cdots, a_i\}}$.
So 
$D^{\bar\bla}=D^{\la^{(1)}}\cdots D^{\la^{(m)}}$ 
is a simple $\sH(\fS_{\bka})$-module, and the set 
$$
\{D^{\bar\bla}\mid \bla\in\La_m^+(r) \text{ $l$-regular},
[\bla]=\bka\} \tag 4.4
$$
forms a complete set of non-isomorphic simple   
$\sH(\fS_{\bka})$-modules.

Since $e_\bka =v_{\bka} T_{w_{\bka'}} c_\bka^{-1}$, where 
 $ T_{w_{\bka'}}$ and $c_\bka$ are invertible and 
$c_\bka$ is in the center of $\sH(\fS_{\bka})$,  \cite{DR, (3.4)} implies that
the map $D^{\bar\bla} \rightarrow
 D^{\bar\bla} e_\bka$ sending $x$ to $x e_\bka$ for 
$x\in D^{\bar\bla}$ is an $R$-module isomorphism. 
Obviously, this is an $\sH(\fS_{\bka})$-module
isomorphism. So $D^{\bar\bla} e_\bka$ is a simple
$\sH(\fS_{\bka})$-module. By (4.3), there is an epimorphism 
from $D^{\bar\bla} e_\bka$ to $D^\bla \epsilon=D^\bla e_\bka $. 
Therefore, 
$D^\bla \epsilon$ is a simple 
$\epsilon \bH \epsilon$-module, and consequently, by (4.4),
all $D^\bla \epsilon$
form a complete set of non-isomorphic $\epsilon \bH \epsilon$-modules. \qed
\enddemo

To state the modular branching rule for Ariki-Koike algebra, we need
the notion of normal and good nodes. Since Ariki-Koike algebra is
semi-simple if $f_{m, r}\neq 0$ and $l>r$ (\cite{A1} or \cite{DR, (5.2)})  
and the branching rule in this semi-simple 
case has been done in (3.9), we assume  in
the rest of this section $l\le r$ and
$f_{m, r}\neq 0$.

Recall from (3.1) that the boxes in the  Young diagram $\sY(\la)$ can be
identified with its position coordinates, called nodes, $\rn=(i, j)$. 
The node $\rn=(i, \la_i)$ is called a removable 
node of $\la$ if $\la_i>\la_{i+1}$ and the node $(i, \la_i+1)$ is called 
an  addable node if $\la_{i-1}>\la_{i}$.
The node $\rn$ is called a removable (resp.  addable) node for
$m$-partition $\bla=(\la^{(1)}, \cdots, \la^{(m)})$
 if it is a removable (resp.  addable) 
node of $\la^{(i)}$ for some $i$.

For each node $(i, j)$ of $\la$, define the $l$-residue of $\rn$
$$\res(\rn)=\cases \text{rem}(j-i) &\text{ if }q=1,\cr
q^{j-i}&\text{ otherwise,}\cr\endcases$$
where $\text{rem}(j-i)$ is the remainder when $j-i$ is divided by $l$.
A removable node $\rn$ of $\la$ is called {\it normal} if for every
addable node $\frak m$ with $\res(\frak m)=\res(\rn)$, there exists
a removable node $\rn'(\frak m)$ strictly between $\rn$ and $\frak m$
with $\res(\rn'(\frak m))=\res(\rn)$, and $\frak m\not=\frak m'$
implies $\rn'(\frak m)=\rn'(\frak m')$.
A removable node is called {\it good} if it is the lowest among
the normal nodes of a fixed residue. (A node $(i,j)$ is
lower than the node $(i',j')$ if $i>i'$.)
Let $\sR_{normal}(\la)$ (resp. $R_{good}(\la)$) be the set of all normal
(resp. good) nodes of $\la\in\La^+(r)$. 
The following result, generalizing Kleshchev's result \cite{K}
for symmetric groups, is  due to Brundan \cite{B, 2.5-6}. 

\proclaim{(4.5) Theorem} Let $\sH(r)$ be the Hecke algebra of type $A$
over the field $R$. If  $\la\in \La^+(r)$  and $\mu\in
\La^+(r-1)$ are  $l$-regular partitions, then

(a)  $\Hom_{\sH(r-1)} (S^\mu, D^\la\vert_{\sH(r-1)})
=\cases R, & \text{ if $\mu=\la_\rn$ for some 
$\rn\in \sR_{normal}(\la)$,}\cr   
       0, & \text{otherwise.}\cr\endcases$     

(b) $\Hom_{\sH(r-1)} (D^\mu, D^\la\vert_{\sH(r-1)})
=\cases R, & \text{ if $\mu=\la_\rn$ for some $\rn\in \sR_{good}(\la)$,}\cr   
       0, & \text{otherwise.}\cr\endcases.$     

Therefore, the socle of the restriction of $D^\la$ to $\sH(r-1)$ is 
$\oplus_{\rn \in  \sR_{good}(\la)} D^{\la_{\rn}}$.  \endproclaim


We are going to generalize this result to the Ariki-Koike algebra
$\bH$ satisfying (4.1).

Let $\bla$ be an $m$-partition of $r$ and $\rn=(i, j)_k$
a node of $\bla$. The residue
of $\rn$ is defined as  
$$
\res(\rn)=\cases\text{rem}(j-i)\xi^k\,\,&\text{ if }q=1, u_k=\xi^k\cr
q^{j-i}u_k &\text{ otherwise,}
\endcases$$
where $\text{rem}(j-i)$ is defined as above
and $\xi$ is the $m$-th primitive root of unity. 
We say that  node  $(i, j)_k$ is {\it lower} than the node 
$(i_1, j_1)_{k_1}$ if either $k>k_1$ or $k=k_1$ and $i>i_1$. 
Thus, any two removable nodes are comparable.
So, we defined normal and good (removable) nodes similarly as in the
$m=1$ case, and let  $\sR_{normal}(\bla)$ (resp. $\sR_{good}(\bla)$) be the set of all normal
(resp. good) nodes of $\bla\in \La_m^+(r)$.

\proclaim{(4.6) Lemma} Let $f_{m, r}\neq 0$ and let $\bla\in \La_m^+(r)$.
If $\rn$ and $\rn'$ are two nodes  of $\sY(\bla)$ with same residues, then 
$\rn, \rn'$ are in $\sY(\la^{(k)})$ for some $k$.
Therefore, $\rn$ is a normal (resp. good) removable node
of $\sY(\bla)$ with  residue
$q^iu_k$ if and only if $\rn$ is a  normal (resp. good) removable node
of $\sY(\la^{(k)})$.\endproclaim

\demo{Proof} Suppose $\res(\rn)=q^au_k$ and $\res(\rn')=q^{a'}u_{k'}$
 with $0\le a,a'\le l-1$ and $1\le k,k'\le m$. Since
$f_{m, r}\neq 0$, $u_i/u_j\neq q^b$ for any $i\neq j$, $1\le i, j\le
m$ and $0\le b\le l-1\le r-1$. Thus,  if $\rn$  and $\rn'$ have the same
residue, that is, $q^au_k=q^{a'}u_{k'}$, then $k=k'$ (and hence, $a=a'$). 
So they are in  
$\sY(\la^{(k)})$.
The last assertion follows immediately from the definition. 
\qed\enddemo

In the proof of modular branching rules below, we shall follow the notations used in \cite{DR}.
Thus, if $\bka=[a_0,a_1,\cdots,a_m]\in\La[m,r]$, then
we define \cite{DR, (1.2)}
$\boldkey a'=[0, r-a_{m-1}, \cdots, r-a_1,r]$, and 
$\boldkey a_\dashv=[a_0, a_1, \cdots, a_{j-1}, r-1, \cdots, r-1],$
where $ j$ is the  minimal index such that
 $a_j=r$. We also define \cite{DR, (1.9)}
 $$
\bka_i=\bka_\dashv+{\boldkey 1}_i,
\text{ where }{\boldkey 1}_i=[0,\undersetbrace{i-1}\to{0,\cdots,0},1,\cdots,1]\in
\La[m,1]. 
$$

\proclaim{(4.7) Lemma}   Let $\bH$ be the Ariki-Koike algebra over 
a commutative ring $R$.
For  $\bka\in \La[m, r]$, write
$\bka_\dashv=[b_0,b_1,\cdots,b_m]\in \La[m, r-1]$, and let 
$ U_i=\bH \pi_{\bka_i}T_{b_i+1, r}
T_{w_{\bka_\dashv}}\tpi_{(\bka_\dashv)'}$ 
 for $1\le i\le m$. Then 

(a) $U_{1}=\bH v_{\bka_1}$ and $U_m =\bH v_{\bka_\dashv}$, and 

(b) $U_i$ is a free $R$-submodule with basis $\{ T_w L_{b_i+1}^c
\pi_{\bka_i} T_{b_i+1, r} T_{w_{\bka_\dashv}}\tpi_{(\bka_\dashv)'}
\mid 0\le c\le i-1, w\in \fS_{r}\}$. So the rank of $U_i$ is $i\cdot
r!$.

 If, in addition, $R$ is an integral domain in which $f_{m, r}$ is a unit, 
then 

(c) $\bH v_{\bka}$ is a projective $\bH$-module, and,
for each $i=2, \cdots, m$, we have  short exact sequence 
$0\rightarrow U_{i-1} \rightarrow
U_i\overset{\ph_i}\to\rightarrow \bH v_{\bka_i}\rightarrow 0$. Therefore,
$\bH v_{\bka_{\dashv}}\cong\oplus_{i=1}^m \bH v_{\bka_i}$.
\endproclaim

\demo{Proof} By the definition, we have
$(\bka_\dashv)'=[0,r-1-b_{m-1},\cdots, r-1-b_1, r-1]$.
So, if $(\bka_\dashv)'=\bkc_\dashv$ for some $\bkc\in\La[m,r]$,
then 
$$\bkc_i=[0,\undersetbrace {i-1}\to {r-1-b_{m-1},\cdots, r-1-b_{m-i+1}},
r-b_{m-i},\cdots,r-b_1,r],$$
and therefore,
$(\bkc_i)'=[0,b_1,\cdots,b_{m-i},b_{m-i+1}+1,\cdots,b_{m-1}+1,r]
=\bka_{m-i+1}$. 
Now, using the $V_i$ defined by $\bkc$ in \cite{DR, (4.7)}, we see that
$$U_i=\Phi(\iota(V_{m-i+1})),\tag4.8$$ 
where $\iota$ is the $R$-linear anti-involution
on $\bH$ sending $T_i$ to $T_i$ (see \cite{GL, (5.5)}) and $\Phi$ is defined
in (2.6). Note that $\Phi(\iota(v_\bkb))=q^av_{\bkb'}$ for some $a\in\BZ$.
Thus, $U_1=\Phi(\iota(V_m))=\Phi(\iota(v_{\bkc_m}\bH))=\bH v_{\bkc_m'}=
\bH v_{\bka_1}$, 
and $U_m=\Phi(\iota(V_1))=\Phi(\iota(v_{\bkc_\dashv}\bH))=
\bH v_{\bkc_\dashv'}=\bH v_{\bka_\dashv}$,
proving (a). (Note that (a) can be also seen directly from the 
definition of $U_i$.)
The statements (b) and (c) follow from \cite{DR, (4.7b),(4.12),(4.14)}
and the relation (4.8).
\qed
\enddemo

Recall from  the proof of (4.2) that we have the isomorphism 
$D^\bla \epsilon \cong D^{\bar \bla}
e_{[\bla]}$ of irreducible $\epsilon\bH\epsilon$-modules. Since
$D^{\bar \bla}e_{[\bla]}= 
(y_{{\la^{(1)}}'\vee\cdots\vee {\la^{(m)}}'} T_{
w(1)^{-1}\cdots w(m)^{-1} } Z_{\bar\bla} v_{[\bla]} \bH) e_{[\bla]}$,
where $Z_{\bar\bla}$ is given in (2.3), 
the  module in parentheses is irreducible, and hence 
isomorphic to $D^\bla$.
In the rest of the section, we put
$$D^\bla= 
y_{{\la^{(1)}}'\vee\cdots\vee {\la^{(m)}}'} T_{
w(1)^{-1}\cdots w(m)^{-1} }Z_{\bar\bla} v_{[\bla]} \bH.\tag 4.9
$$ 
Recall also that $\La_m^+(r)$ is the set of all $m$-partitions of $r$. The
following result is the modular branching rule for the Ariki-Koike
algebra under the assumption (4.1), i.e., $f_{m, r}\neq 0$ in $R$.

\proclaim{(4.10) Theorem} Let $\bH_{m}^r$ be the Ariki-Koike algebra
over a field $R$ in which $f_{m, r}\neq 0$.  
 If  $\bla\in \La_m^+(r)$  and $\brho\in
\La_m^+(r-1)$ are  $l$-regular $m$-partitions, then

(a)  $\Hom_{\bH_{m}^{r-1}} (S^\brho, D^\bla\vert_{\bH_{m}^{r-1}})
\cong \cases R, & \text{ if $\brho=\bla_\rn$ for some $\rn\in
\sR_{normal}(\bla)$,}\cr   
       0, & \text{ otherwise.}\cr\endcases$     

(b) $\Hom_{\bH_{m}^{r-1}} (D^\brho, D^\bla\vert_{\bH_{m}^{r-1}})
\cong \cases R, & \text{ if $\brho=\bla_\rn$ for some $\rn\in
\sR_{good}(\bla)$,}\cr   
       0, & \text{otherwise.}\cr\endcases.$     

Therefore, the socle of the restriction of $D^\bla$ to $\bH_{m}^{r-1}$ is 
$\oplus_{\rn \in  \sR_{good}(\bla)} D^{\bla_{\rn}}$.  \endproclaim

\demo{Proof}  
Because $f_{m, r-1}$ is a factor of $f_{m, r}$ and 
$f_{m, r}\neq 0$, we have  
$f_{m-1, r}\neq 0$. So, by \cite{DR, (4.14c)}, we have
a Morita equivalence between the categories of $\bH_{m}^{r-1}$-modules
and $\epsilon \bH_{m}^{r-1}\epsilon$-modules, where 
$\epsilon=\sum_{\bkb\in \La[m, r-1]} e_\bkb$. Note that 
$\{e_\bkb\}$ is a set of orthogonal idempotents. 
Using a standard result on Morita equivalence (see, e.g., \cite{AF})
and noting \cite{DR, (3.10)},
we have, for  $M=S^\brho$ or $D^\brho$,  the linear isomorphism 
$$\aligned 
\Hom_{\bH_{m}^{r-1}} (M, D^\bla\vert_{\bH_{m}^{r-1}}) & \cong    \Hom_{\epsilon  \bH_{m}^{r-1}\epsilon }
(M\epsilon   , D^\bla\epsilon)\cr & \cong
\Hom_{e_{[\brho]}  \bH_{m}^{r-1} e_{[\brho]} }
(M e_{[\brho]},  D^\bla e_{[\brho]}).\cr\endaligned\tag4.11
$$
Since 
$e_{[\brho]}  \bH_{m}^{r-1} e_{[\brho]}=e_{[\brho]}\sH(\fS_{[\brho]})=
\sH(\fS_{[\brho]})e_{[\brho]}
\cong \sH(\fS_{[\brho]})$
we may twist the action on $M e_{[\brho]}$ via this isomorphism and obtain
isomorphisms of $\sH(\fS_{[\brho]})$-modules: 
$S^\brho e_{[\brho]}=S^{\bar\brho} e_{[\brho]}\cong S^{\bar\brho}$
and
$D^\brho e_{[\brho]}=D^{\bar\brho} e_{[\brho]}\cong D^{\bar\brho}$.
So 
the calculation
of the Hom sets above is reduced 
to calculate the module $D^\bla e_{[\brho]}$. 

From (4.9), we
first have to manipulate the set $v_{[\bla]}\bH e_{[\brho]}
=v_{[\bla]}\bH v_{[\brho]}T_{w_{\bka'}}c_{[\brho]}^{-1}$
(see line 2 after (4.4)). Choose $\bka\in \La[m, r]$ such that 
$\bka_\dashv=[\brho]$. Then $v_{[\bla]}\bH v_{[\brho]}=v_{[\bla]}\bH v_{\bka_\dashv}$. Suppose $v_{[\bla]}\bH v_{\bka_\dashv}\neq0$. Then, the isomorphism
$\bH v_{\bka_\dashv}\cong  \oplus_{i=1}^m \bH v_{\bka_i}$ in (4.7c) together
with the orthogonal property  \cite{DR, (3.10)}
implies that there is a unique $i$ such that $[\bla]=\bka_i$.
Thus, $v_{[\bla]}\bH v_{\bka_\dashv}=v_{\bka_i}\bH v_{\bka_\dashv}=
v_{\bka_i}U_i$, since $U_i\cong \oplus_{j=i}^m\bH v_{\bka_j}$.
Now $\tilde\pi_{\bka_i'}\sH\pi_{\bka_i}=\sH(\fS_{\bka_i'})\iota(v_{\bka_i})$
by \cite{DR, (3.1)}, it follows that  
$$
\aligned
v_{\bka_i}U_i&=\pi_{\bka_i}T_{w_{\bka_i}}(\tilde\pi_{\bka_i'}\bH
\pi_{\bka_i})T_{b_i+1, r}T_{w_{\bka_\dashv}}\tpi_{(\bka_\dashv)'}\cr
&=\sH(\fS_{\bka_i})\pi_{\bka_i}
T_{w_{\bka_i}}\iota(v_{\bka_i})T_{b_i+1,r}T_{w_{\bka_\dashv}}
\tilde\pi_{(\bka_\dashv)'}\cr
&=\sH(\fS_{\bka_i})v_{\bka_i}
T_{w_{\bka_i'}}\pi_{\bka_i}T_{b_i+1,r}T_{w_{\bka_\dashv}}
\tilde\pi_{(\bka_\dashv)'}.\cr
\endaligned$$
Write  $\pi_{\bka_i}T_{b_i+1,r}=h_i\pi_{\bka_{\dashv}}$, where
$$
\aligned
h_i&=(L_{b_{i}+1}-u_{i+1}) 
T_{b_{i}+1, b_{i+1}+1}\cdots (L_{b_{m-1}+1}-u_m) T_{b_{m-1}+1, r}\cr
&=\prod_{j=i}^{m-1}(L_{b_{j}+1}-u_{j+1}) 
T_{b_{j}+1, b_{j+1}+1}\cr
\endaligned
$$
(see \cite{DR, (2.7b)}).
Then, $v_{\bka_i}\bH v_{\bka_\dashv}=\sH(\fS_{\bka_i})v_{\bka_i}T_{w_{\bka_i'}}h_iv_{\bka_{\dashv}}$, and therefore, $D^\bla e_{\bka_\dashv}=D^{\bar\la}
v_{\bka_i}T_{w_{\bka_i'}}h_ie_{\bka_{\dashv}}$.
Thus the isomorphism $\sH(\fS_{\bka_\dashv})\cong 
e_{\bka_{\dashv}}\bH^{r-1}_m e_{\bka_\dashv}=e_{\bka_{\dashv}}
\sH(\fS_{\bka_\dashv})$ 
 will turn as above  
the $e_{\bka_{\dashv}}
\sH(\fS_{\bka_\dashv})$-module
$D^{\bar\la}
v_{\bka_i}T_{w_{\bka_i'}}h_ie_{\bka_{\dashv}}$ into an 
$\sH(\fS_{\bka_\dashv})$-module with the action
$$(xv_{\bka_i}T_{w_{\bka_i'}}h_ie_{\bka_{\dashv}})*T_w=xv_{\bka_i}T_{w_{\bka_i}'}h_iT_we_{\bka_{\dashv}}\tag 4.12$$ 
for all $x\in D^{\bar\la}$ and $w\in \fS_{\bka_\dashv}$.

We now claim that this $\sH(\fS_{\bka_\dashv})$-module is isomorphic to the
$\sH(\fS_{\bka_\dashv})$-module $D^{\bar\la} T_{b_i+1, r}$.
Indeed, write $w=w_1\cdots w_m$, where $w_j\in\fS_{\{b_j+1,\cdots,b_{j+1}\}}$.
Since 
$$s_{b_j+1,b_{j+1}+1}\fS_{\{b_j+1,\cdots,b_{j+1}\}}=\fS_{\{b_j+2,\cdots,
b_{j+1}+1\}}s_{b_j+1,b_{j+1}+1}$$ 
and the product has length additivity,
we have 
$$T_{b_j+1,b_{j+1}+1}T_{w_j}=T_{w_j^*}T_{b_j+1,b_{j+1}+1}$$
where $w_j^*=s_{b_j+1,b_{j+1}+1}w_js_{b_{j+1}+1,b_{j}+1}$.
Thus, from (4.12), we have by the commuting relations between $L_i$'s and $T_j$'s
(\cite{DR, (2.5)})
$$\aligned
xv_{\bka_i}T_{w_{\bka_i'}}h_iT_we_{\bka_{\dashv}}
&=
xv_{\bka_i}T_{w_{\bka_i'}}T_{w_1,\cdots w_{i-1}}\prod_{j=i}^{m-1}(L_{b_{j}+1}-u_{j+1}) (T_{b_{j}+1, b_{j+1}+1}T_{w_j})\cr
&=xv_{\bka_i}T_{w_{\bka_i'}}T_{w_1,\cdots w_{i-1}w_i^*\cdots w_m^*}h_i\cr
&=xT_{w_1,\cdots w_{i-1}w_i^*\cdots w_m^*}
v_{\bka_i}T_{w_{\bka_i}'}h_ie_{\bka_{\dashv}}.\cr
\endaligned
$$
Here the last equality follows from a repeated use of \cite{DR, (1.8)}
and \cite{DR, (2.7a)} noting $w_1,\cdots w_{i-1}w_i^*\cdots w_m^*\in\fS_{\bka_i}$. On the other hand, we have
$$xT_{b_i+1,r}T_w=xT_{w_1,\cdots w_{i-1}w_i^*\cdots w_m^*}T_{b_i+1,r}.\tag4.13$$
Therefore, the linear map 
$$f:D^{\bar\la} T_{b_i+1, r}\to D^{\bar\la}
v_{\bka_i}T_{w_{\bka_i'}}h_ie_{\bka_{\dashv}}$$
defined by sending $xT_{b_i+1,r}$ to $xv_{\bka_i}T_{w_{\bka_i}'}h_ie_{\bka_{\dashv}}$ is a surjective  $\sH(\fS_{\bka_\dashv})$-module homomorphism,
and hence an isomorphism by a comparison of dimensions, proving the claim.

Let  
$[\la]=\bka_i=[c_0,\cdots,c_m]$. Then,
$c_j=b_j$ if $0\le j\le i-1$, or $b_j+1$ if $i\le j\le m$.
Since $D^{\bar\la}=D^{\la^{(1)}}\cdots D^{\la^{(m)}}$ where
$D^{\la^{(j)}}$ is the corresponding irreducible 
$\sH(\fS_{\{c_{j-1}+1,\cdots,c_{j}\}})$-module,
(4.13) implies the $\sH(\bka_\dashv)$-module isomorphism 
$$
D^{\bar\la} T_{b_i+1, r}\cong D^{\la^{(1)}}\otimes\cdots\otimes D^{\la^{(i-1)}}
\otimes(D^{\la^{(i)}}\vert_{\sH(\fS_{b_{i-1}+1,\cdots,b_{i}})})\otimes
\tilde D^{\la^{(i+1)}}\otimes\cdots\otimes \tilde D^{\la^{(r)}},
$$
where  $D^{\la^{(i)}}\vert_{\sH(\fS_{\{b_{i-1}+1,\cdots,b_{i}})\}}$
denotes the module by restricting the action from $\sH(\fS_{\{c_{i-1}+1,\cdots,c_{i}\}})$ to
$\sH(\fS_{c_{i-1}+1,\cdots,b_{i}})$, and
$\tilde D^{\la^{(j)}}$, for $j=i+1,\cdots,m$, denotes the irreducible
$\sH(\fS_{\{b_{j-1}+1,\cdots,b_{j}\}})$-module obtained by twisting the
$\sH(\fS_{\{c_{j-1}+1,\cdots,c_{j}\}})$-module structure on $D^{\la^{(j)}}$
through the canonical isomorphism $\sH(\fS_{\{c_{j-1}+1,\cdots,c_{j}\}})\cong
\sH(\fS_{\{b_{j-1}+1,\cdots,b_{j}\}})$.

With what we obtained above, (4.11) becomes (recalling $[\bla]=\bka_i$)
$$
\aligned & \Hom_{e_{[\brho]} \bH_{m}^{r-1} e_{[\brho]}} 
(S^{\bar \brho} e_{[\brho]},  D^{\bla} e_{[\brho]})
\cong \Hom_{ \sH (\fS_{[\brho)]})  } (S^{\bar \brho}, D^{\bar\bla}T_{b_i+1,r} )\cr  &
\cong \otimes_{j=1}^{i-1} 
\Hom_{\sH(\fS_{\{b_{j-1}+1,\cdots, b_j\} })}( S^{\mu^{(j)}},
D^{\la^{(j)}} ) \cr &
\quad\,\otimes \Hom_{\sH(\fS_{\{b_{i-1}+1,\cdots, b_i \}})}( S^{\mu^{(i)}},
D^{\la^{(i)}}\vert_{\sH(\fS_{\{b_{i-1}+1,\cdots,b_{i}\}})})\otimes \cr
& \quad\,\otimes_{j=i+1}^m 
\Hom_{\sH(\fS_{\{b_{j-1}+1,\cdots, b_j\} })} ( S^{\mu^{(j)}},
\tilde D^{\la^{(j)}}),\cr\endaligned\tag 4.14
$$
and similarly,
$$
\aligned & \Hom_{e_{[\brho]} \bH_{m}^{r-1} e_{[\brho]}} 
(D^{\bar \brho} e_{[\brho]},  D^{\bla} e_{[\brho]})
\cong \Hom_{ \sH (\fS_{[\brho]})  } (D^{\bar \brho}, D^{\bar\bla} T_{b_i+1,
r})\cr  &
\cong \otimes_{j=1}^{i-1} 
\Hom_{\sH(\fS_{\{b_{j-1}+1,\cdots, b_j\} })}( D^{\mu^{(j)}},
D^{\la^{(j)}} ) \cr &
\quad\,\otimes \Hom_{\sH(\fS_{\{b_{i-1}+1,\cdots, b_i\}})}( D^{\mu^{(i)}},
D^{\la^{(i)}}\vert_{\sH(\fS_{\{b_{i-1}+1,\cdots,b_{i}\}})})\otimes \cr
&\quad\, \otimes_{j=i+1}^m 
\Hom_{\sH(\fS_{\{b_{j-1}+1,\cdots, b_j\} }) } ( D^{\mu^{(j)}},
\tilde D^{\la^{(j)}}).\cr\endaligned \tag 4.15
$$
Clearly, 
the Hom set in (4.14) (resp. (4.15)) is isomorphic to $R$ if and only if 
$\mu^{(j)}=\la^{(j)}$ for all $j\neq i$, and $\Hom_{\sH(\fS_{\{b_{i-1}+1,\cdots, b_i\} })}( S^{\mu^{(i)}},
D^{\la^{(i)}}\vert_{\sH(\fS_{\{b_{i-1}+1,\cdots,b_{i}\}})} )\cong R$
(resp. $\Hom_{\sH(\fS_{\{b_{i-1}+1,\cdots, b_i\} })}( D^{\mu^{(i)}},
D^{\la^{(i)}}\vert_{\sH(\fS_{\{b_{i-1}+1,\cdots,b_{i}\}})} )\cong R$).
However, by (4.5), the latter is equivalent to 
$\mu^{(i)}={\la_\rn^{(i)}}$ for 
$\rn\in
\sR_{normal}(\la^{(i)})$ (resp. $\rn\in
\sR_{good}(\la^{(i)})$). Now, applying (4.6), the theorem is proved.\qed
\enddemo 

Using the very recent work \cite{A2}
on the classification of the irreducible modules, we propose the 
following conjecture.

\proclaim{(4.16) Conjecture} Replacing the $l$-regular multipartitions
by Kleshchev multipartitions, (4.10) holds in general.\endproclaim

\end

From sfb8@Mathematik.Uni-Bielefeld.DE Thu Dec 17 03:55:30 1998
Received: from crash.Mathematik.Uni-Bielefeld.DE (crash.Mathematik.Uni-Bielefeld.DE [129.70.14.10])
        by alpha.maths.unsw.EDU.AU (8.8.8/8.8.8) with ESMTP id DAA14582
        for <jied@maths.unsw.edu.au>; Thu, 17 Dec 1998 03:55:22 +1100 (EST)
Received: from charybdis.Mathematik.Uni-Bielefeld.DE (charybdis.Mathematik.Uni-Bielefeld.DE [129.70.15.10])
        by crash.Mathematik.Uni-Bielefeld.DE (8.8.8/8.8.8) with ESMTP id RAA02667;
        Wed, 16 Dec 1998 17:54:35 +0100 (MET)
From: SFB <sfb8@Mathematik.Uni-Bielefeld.DE>
Received: (from sfb8@localhost)
        by charybdis.Mathematik.Uni-Bielefeld.DE (8.8.3/8.8.8) id RAA04462;
        Wed, 16 Dec 1998 17:54:35 +0100 (MET)
Date: Wed, 16 Dec 1998 17:54:35 +0100 (MET)
Message-Id: <199812161654.RAA04462@charybdis.Mathematik.Uni-Bielefeld.DE>
To: jied@maths.unsw.edu.au
Subject: tex file
Cc: sfb8@Mathematik.Uni-Bielefeld.DE
Status: O

I add some words to illustarte the display (4.12). I think it is better to
use the new version to replace the old old. Similarly, I rewrite the 
display (4.13). Now one may see easily how to prove (a), just the same as (b).

I will go back to China on December 30. My wife told me that the ECNU sent
a document to our university about my moving. I will move to ECNU in the early
of 1999 if everything is ok. Do you have time to visit ECNU this year? 

I can make an arrangement since I have got the NNSF for 1999-2001.
Best, Hebing

The proof of (d) is similar to the one given in \cite{DR, (4.12-3)}.
In fact, 
the proof below shows that the restriction $i\ge k$ there can be removed.
 First, we define a homomorphism
$f_i:U_i\to \bH v_{\bka_i}$. For this purpose, we write by \cite{DR, 1.10}
$$
v_{\bka_i}
=\pi_{\bka_i} T_{w_{\bka_i}} \tpi_{\bka_i'}
=\pi_{\bka_i} T_{b_i+1} T_{w_{\bka_\dashv} }
(T_{r, r-b_{i-1}} \tpi_{\bka_i'}).
$$
Since $\bka_i=[\,\cdots,b_{i-1}, b_i+1, \cdots, b_m+1]$ (here $b_m=r-1$), we have 
$\bka_i'=[0, b_m-b_{m-1}, \cdots, b_m-b_i, r-b_{i-1}, \cdots, r-b_1, r]$
and $\bka_\dashv'=[0, b_m-b_{m-1}, \cdots, b_m-b_1, b_{m}]$, and by 
\cite{DR, (2.5-6)}, 
$T_{r, r-b_{i-1} } \tpi_{\bka_i'}=\tpi_{(\bka_\dashv)'} h_i$, where
$$h_i=\prod_{l=1}^{k-1} (L_r-u_l) \times
T_{r, r-b_k} (L_{r-b_k}-u_k)\cdots T_{r-b_{i-2}, r-b_{i-1}} 
(L_{r-b_{i-1}}-u_{i-1}).$$
Here $k$ is the minimal number with $b_k\neq 0$. Thus, we define
$f_i$ by $f_i(x)=xh_i$ for all $x\in U_i$. It is clearly
a homomorphism of (left) $\bH$-modules.

{\bf Complete the proof! In the proof of [DR, 4.13], condition
$i\ge k$ is used. How can you avoid this here. YOU MUST SHOW THIS.}
\qed\enddemo